\documentclass[12pt]{article}

\usepackage{latexsym}
\newtheorem{teo}{Theorem}[section]
\newtheorem{ex}[teo]{Example}
\newtheorem{definicao}[teo]{Definition}
\newtheorem{lema}[teo]{Lemma}
\newtheorem{prop}[teo]{Proposition}
\newtheorem{cor}[teo]{Corollary}
\newtheorem{obs}[teo]{Remark}

\newcommand{\grad}{\ensuremath{\mathrm{grad}\ }}

\newcommand{\tub}{\ensuremath{\mathrm{Tub} }}

\newcommand{\F}{\ensuremath{\mathcal{F}}}
\newcommand{\singularF}{\ensuremath{\mathcal{X}_{F}}}
\newcommand{\rank}{\ensuremath{\mathrm{rank}\ }}

\newcommand{\agradecimentos}{\textsf{Acknowledgmets:\ }} 
\newcommand{\dem}{\textbf{Proof:\  }}
\newcommand{\afir}{\textsf{Statement \  }}

\newcommand{\fimdem}{\Box}

\newcommand{\Holsing}{\ensuremath{\mathrm{Holsing}}}
\newcommand{\Iso}{\ensuremath{\mathrm{Iso}}}

\begin{document}

%first part of the article: local and global aspects

\title{Singular Riemannian Foliations with Sections \footnotemark\footnotetext{2000 Mathematics Subject Classifications. 53C12, 57R30\\Key words and 
phrases. Singular riemannian foliations, isoparametric maps, equifocal submanifolds, isoparametric submanifolds, singular holonomy.}
 }
\author{Marcos M. Alexandrino\thanks{Supported by CNPq 
}\\ 
\small Departamento de Matem\'{a}tica,\\
\small Pontif\'{\i}cia Universidade Cat\'{o}lica,\\ 
\small Rua Marqu\^{e}s de S\~{a}o Vicente, 225\\
\small  22453-900, Rio de Janeiro, Brazil.\\
\small email: malex@mat.puc-rio.br
}
\date{November 2003}

\maketitle

\begin{abstract}
A singular foliation on a complete riemannian manifold is said to be riemannian if  every geodesic that is perpendicular at one point to a leaf remains perpendicular to every leaf it meets. 
In this paper we study singular riemannian foliations that have sections, i.e., totally geodesic complete immersed submanifolds that meet each leaf orthogonally and whose dimensions are the codimensions of the regular leaves.

We prove here  that the 
restriction of the foliation to a slice of a leaf is diffeomorphic  to an isoparametric foliation on an open set of  an euclidian space.  
This result gives us local information about the singular foliation and in particular  about the singular stratification of the foliation. It also allows us to  describe the plaques of  the foliation as level sets of a transnormal map (a generalisation of an isoparametric map).
We also prove that  the regular leaves of a singular riemannian foliation with sections  are locally equifocal. We use this property to define a singular holonomy. Then we  establish some results about this singular holonomy and illustrate them with a couple of  examples.
 
\end{abstract}

\section{Introduction}
In this section we shall introduce the concept of a singular riemannian foliation with sections, review typical examples of this kind of foliation
and  state our main results (Theorem \ref{frss-eh-equifocal} and Theorem \ref{sliceteorema}), which relate the new concept with the concepts of isoparametric and equifocal submanifolds.

We start by recalling  the   definition of a singular riemannian foliation (see  the  book of  P. Molino \cite{Molino}).

\begin{definicao}
\emph{ A partition $\F$ of a complete riemannian manifold $M$ by connected immersed submanifolds (the \emph{leaves}) is called \emph{singular riemannian foliation on $M$} if it verifies the following conditions
\begin{enumerate}
\item $\F$ is \emph{singular}, i.e., the set $\singularF$ of smooth vector fields on $M$  that are tangent at each point to the corresponding leaf is transitive on each leaf. In other words, for each leaf $L$ and each $p\in L,$ one can find vector fields $v_{i}\in \singularF$ such that $\{v_{i}(p)\}$  is a basis of $T_{p}L.$
\item  The partition is \emph{transnormal}, i.e., every geodesic that is perpendicular at one point to a leaf remains perpendicular to every leaf it meets.
\end{enumerate}
}
\end{definicao}

Let $\F$ be a singular riemannian foliation on an complete riemannian manifold $M.$  A point $p\in M$ is called \emph{regular} if the dimension of the leaf $ L_{p}$ that contains $p$ is maximal.  A point is called \emph{singular} if it is not  regular.
Let $L$ be an immersed submanifold of a riemannian manifold $M.$  A section $\xi$ of  the normal bundle $\nu (L)$ is said  to be a \emph{parallel normal field} along $L$ if  $\nabla^{\nu}\xi\equiv 0,$ where $\nabla^{\nu}$ is the normal connection.  $L$  is said to have globally flat normal bundle,  if the holonomy  of the normal bundle  $\nu (L)$ is trivial, i.e., if any normal vector can be extended to a globally defined parallel normal field.

\begin{definicao}[s.r.f.s.]
\emph{ Let $\F$ be a singular riemannian foliation on a complete riemannian manifold $M.$
$\F$ is said to be a \emph{singular riemannian foliation with section} (s.r.f.s. for short) if for every regular point $p,$ the set $\sigma :=\exp_{p}(\nu L_{p})$ is an immersed  complete submanifold that meets each leaf orthogonally and if the regular points of $\sigma$ are dense in it.  $\sigma$ is called a \emph{section}.}
\end{definicao}

Let $p\in M$ and $\tub(P_{p})$ be a tubular neighborhood of a plaque $P_{p}$ that contains $p.$ Then the connected component of $\exp_{p}(\nu P_{p})\cap \tub(P_{p})$  that contains $p$ is called a \emph{slice} at $p.$ Let $\Sigma_{p}$ denote it. 
Now consider the intersection of $\tub(P_{p})$ with a section of the foliation. Each connected component of this set is called a \emph{local section}. 
These two concepts play here an important role and are  related to each other. In fact,  we show in  Proposition \ref{slice-uniao-secoes} that \emph{ the slice at a singular point is the union of the local sections that contain this singular point.}

Typical examples of singular riemannian foliations with sections are the orbits of a polar action, parallel submanifolds of an isoparametric submanifolds in a space form and parallel submanifolds of an equifocal submanifold with flat sections in a compact symmetric space, concepts that we now recall.

An isometric action of a compact Lie group on a riemannian manifold $M$ is called \emph{polar} if there exists a complete immersed submanifold $\sigma$ of $M$ that meets each $G$-orbit orthogonally. Such $\sigma$ is called a  \emph{section}. A typical example of a polar  action is a compact Lie group with a biinvariant metric that acts on itself  by conjugation. In this case the maximal tori are the sections.

A submanifold of a real space form is called \emph{isoparametric} if its normal bundle is flat and if the principal curvatures along any parallel normal vector field are constant. The history of isoparametric hypersurfaces  and submanifolds and their generalizations can be found in the survey \cite{Th} of G. Thorbergsson. 

Now we recall the concept of an equifocal submanifold that was introduced by C.L. Terng and G. Thorbergsson \cite{TTh1} as a generalization of the concept of an isoparametric submanifold.

\begin{definicao}
\emph{
A connected immersed submanifold  $L$  of a complete riemannian manifold $M$ is called \emph{equifocal} if
\begin{enumerate}
\item[0)] the normal bundle $\nu(L)$ is globally flat,
\item[1)] for each  parallel normal field $\xi$  along $L,$  the derivative of  the map    $\eta_{\xi}:L\rightarrow M,$ defined as $\eta_{\xi}(x):=\exp_{x}(\xi),$ has constant rank,
\item[2)] $L$ has sections, i.e., for all $p\in L$ there exists a complete, immersed, totally geodesic submanifold $\sigma$  such that $\nu_{p}(L)=T_{p}\sigma.$
\end{enumerate}
}
\end{definicao}
A connected immersed submanifold $L$ is called \emph{locally equifocal} if, for each $p\in L,$ there exists a neighborhood $U\subset L$ of $p$ in $L$ such that $U$ is an equifocal submanifold.
%%%%%%%%pre-introducao

Finally we are ready to state our main results.

\textbf{Theorem \ref{frss-eh-equifocal}}
\emph{The regular leaves of a singular riemannian foliation with sections on a complete riemannian manifold $M$ are locally equifocal. In addition, if  
all the leaves are compact, then 
the union of  regular leaves that are equifocal is an open and dense set in $M$.}

\vspace{4mm}

This result implies that given an equifocal leaf $L$ we can reconstruct the singular foliation taking all parallel submanifolds of $L$ (see Corollary \ref{cor-map-paralelo}). In other words, \emph{let $L$ be a regular equifocal leaf and $\Xi$ denote the set of all parallel normal fields along $L.$ Then $\F=\{\eta_{\xi}(L)\}_{\xi\in \, \Xi}.$} 
Theorem \ref{frss-eh-equifocal} allows us to define a singular holonomy. We also establish some results about this singular holonomy (see section 3) and illustrate them with a couple of new examples. Theorem \ref{frss-eh-equifocal} is also used to prove the following result:

\vspace{4mm}
\textbf{Theorem \ref{sliceteorema} (slice theorem)} 
\emph{ Let $\F$ be a singular riemannian foliation with sections  on a complete riemannian manifold $M$  and  $\Sigma_{q}$ the slice at a point  $q\in M.$ Then  $\F$ restricted to  $\Sigma_{q}$ is diffeomorphic to an isoparametric foliation on an open set of  $\mathbf{R}^{n},$  where $n$ is the dimension of   $\Sigma_{q}.$ }
\vspace{4mm}

Owing to the slice theorem, we can see the plaques of the singular foliation, which are in a tubular neighborhood of a singular plaque $P,$ as the product of isoparametric submanifolds  and  $P.$ In particular, we can better understand the singular stratification (see Corollary \ref{estratificacao-singular}).

%%%%%%%%%%%%principais resultados%%
A consequence of the slice theorem is  Proposition  \ref{frss-transnormal}  that claims that   \emph{ the plaques of a s.r.f.s. are always level sets of a transnormal map,} concept that we recall below.

\begin{definicao}
\emph{
Let  $M^{n+q}$  be a complete riemannian manifold. A smooth map   $H=(h_{1}\cdots h_{q}):M^{n+q}\rightarrow \mathbf{R}^{q}$ is called  \emph{transnormal} if
\begin{enumerate}
\item[0)]  $H$ has a regular value,
\item[1)]  for every regular value $c$ there exists a neighborhood $V$ of $H^{-1}(c)$ in $M$ and smooth functions $b_{i\,j}$ on $H(V)$ such that, for every $x\in V,$ \linebreak
$<\grad h_{i}(x),\grad h_{j}(x)>= b_{i\,j}\circ H(x),$
\item[2)]  there is a sufficiently small neighborhood of each  regular level set such that $[\grad h_{i},\grad h_{j}]$ is a linear combination of $\grad h_{1}\cdots \grad h_{q},$ with coefficients being functions of $H,$ for all $i$ and $j$. 
\end{enumerate}
}
\end{definicao}
This definition is equivalent to saying that $H$ has a regular value and for each regular value $c$ there exists a neighborhood  $V$ of $ H^{-1}(c)$ in $M$ such that $H \mid_{V}\rightarrow H(V)$ is an integrable riemannian submersion, where the metric $(g_{i\,j})$ of  $H(V)$ is the inverse matrix of $(b_{i\,j}).$

A transnormal map $H$ is said to be an \emph{isoparametric map} if $V$ can be chosen to be $M$  and   $\bigtriangleup h_{i}= a_{i}\circ H,$ where $a_{i}$ are smooth functions. 
 
Isoparametric submanifolds in  space forms and equifocal submanifolds with flat sections in  simply connected symmetric spaces of compact type can always be  described as regular level sets of transnormal analytic maps,  see R.Palais and C.L.Terng \cite{PTlivro} and E. Heintzte, X.Liu and C.Olmos \cite{HOL}. 

We prove in \cite{Ale} that \emph{the regular leaves of an analytic transnormal map on an analytic complete manifold are equifocal submanifolds and leaves of a singular riemannian foliation with sections}. Hence,  Proposition \ref{frss-transnormal} is a local converse of this result.

This paper is organized as follows. In section 2 we shall prove some propositions about singular riemannian foliation with sections (s.r.f.s. for short), Theorem \ref{frss-eh-equifocal} and Theorem \ref{sliceteorema}. In section 3 we shall introduce the concept of singular holonomy of a s.r.f.s. and establish some results about it. In section 4 we illustrate some properties of  singular holonomies constructing singular foliations by suspensions of  homomorphisms. 

\agradecimentos This paper is part of my PhD thesis \cite{AlePhD}. I would like to thank my thesis advisors Prof. Ricardo Sa Earp (PUC-Rio) and Prof. Gudlaugur Thorbergsson (Uni zu K\"{o}ln) for their consistent support during my study in Brazil and Germany and for many helpful discussions. I also thank my friend Dirk T\"{o}ben for useful suggestions.

%%%%%%%%fim
%%%%%%%%%%%%%%%%%%%%%%%
%%%%%%%%%%%%%%%RESULTS
\section{Proof of the main results}

%%%%%%%%%%%%%%%PROP: SLICE EH A UNIAO DE SECOES%%%%%%%

\begin{prop}
\label{slice-uniao-secoes}
Let $\F$ be a s.r.f.s. on a complete riemannian manifold $M$ and let  $q\in M.$  Then 
\begin{enumerate}
\item[a)] $\Sigma_{q}= \cup_{\sigma\in\Lambda (q)}\, \sigma,$ where $\Lambda(q)$  is the set of all local sections that contain $q.$
\item[b)] $\Sigma_{x}\subset\Sigma_{q}$ for all $x\in\Sigma_{q}.$
\item[c)]  $T_{x}\Sigma_{q}=T_{x}\Sigma_{x}\oplus T_{x}(U\cap \Sigma_{q}),$ where  $x\in\Sigma_{q}$ and  $U\subset L_{x}$ is an open set of  $x$ in $L_{x}.$ 
\end{enumerate}
\end{prop} 

\dem a)  At first we check that $\Sigma_{q}\supset \cup_{\sigma\in\Lambda (q)}\, \sigma.$ Let $\sigma$ be a local section that contains $q,$  let $p$ be a regular point of $\sigma$  and $\gamma$ the shortest segment of geodesic   that joins $q$ to $p.$ Then  $\gamma$ is orthogonal to $L_{p}$ for $\gamma\subset\sigma$ and $\sigma$ is orthogonal to $L_{p}.$ Since $\F$ is a riemannian foliation, $\gamma$ is also orthogonal to $L_{q}$ and hence $p\in\Sigma_{q}.$ Since the regular points are dense in $\sigma,$ $\Sigma_{q}\supset\sigma.$

Now we check that $\Sigma_{q}\subset  \cup_{\sigma\in\Lambda (q)}\, \sigma.$  Let $p\in\Sigma_{q}$ be a regular point and $\gamma$ the  segment of geodesic orthogonal to $L_{q}$   that joins $q$ to $p.$ Since $\F$ is a riemannian foliation, $\gamma$ is orthogonal to $L_{p}.$ Therefore $\gamma$ belongs to the local section $\sigma$ that contains $p.$ In particular $q\in \sigma.$  In other words, each regular point $p\in\Sigma_{q}$ belongs to a local section $\sigma$ that contains $q.$
 
Finally let $z\in\Sigma_{q}$ be a singular point, $\sigma$ a local section that contains $z,$ and $p$ a regular point of $\sigma.$ Since the slice is defined on a tubular neighborhood of a plaque $P_{q},$ there exists only one point $\tilde{q}\in P_{q}$ such that $p\in\Sigma_{\tilde{q}}.$ As we have shown above, $\tilde{q}\in\sigma.$  Now it follows  from the first part of the proof that   $z\in\sigma\subset\Sigma_{\tilde{q}}$ Since $z\in\Sigma_{q},$   $\tilde{q}=q.$  

b) Let $x\in\Sigma_{q}$ and $\sigma\subset\Sigma_{x}$ a local section. It follows from the proof of item a) that $\sigma\subset \Sigma_{q}$ and $q\in\sigma.$ Since $\Sigma_{x}$ is a union of local sections that contain $x,$ $\Sigma_{x}\subset\Sigma_{q}.$

c) Since the foliation $\F$ is singular, we have:

\begin{eqnarray*}
\dim T_{x}\Sigma_{x}+ \dim T_{x}(U\cap \Sigma_{q})+\dim L_{q}&=&\dim M\\
                                                             &=&\dim T_{x}\Sigma_{q}+ \dim L_{q}.
\end{eqnarray*}
The item b)  and the above equation imply the item c) $\fimdem$

%%%%%%%%%%%%%%%%%%%%COROLARIO: FOLHEACAO NO SLICE EH FRSS
\begin{obs}
\emph{
In \cite{Molino} (pag 209) Molino showed that given a singular riemannian foliation  $\F$ it is possible to change the metric in such way that the restriction of $\F$ to a slice is a singular riemannian foliation with respect to this new metric. This change respect the distance between the leaves. As we see below this change is not necessary if the singular riemannian foliation has sections.  
}
\end{obs}

\begin{cor}
Let  $\F$ be a s.r.f.s. on a complete riemannian manifold  $M$ and  $\Sigma$  a slice. Then  $\F\cap\Sigma$  is a s.r.f.s. on $\Sigma$ with the induced metric of $M.$
\end{cor}
\dem Let $\gamma$ be a segment of geodesic that is orthogonal to $L_{x}$  where $x\in\Sigma.$ Since $\gamma\subset \Sigma_{x},$ it follows from the item b) of Proposition \ref{slice-uniao-secoes}
 that  $\gamma\subset \Sigma.$ Since $\F$ is riemannian, $\gamma$ is orthogonal to the leaves of $\F\cap\Sigma.$ Therefore $\F\cap\Sigma$ is a singular riemannian foliation.  
  
Now let $\sigma$ be a local section that contains $x.$ Then it follows from item a) of Proposition \ref{slice-uniao-secoes}
 that $\sigma\subset \Sigma$ and hence $\sigma$ is a local section of $\F\cap\Sigma$. Therefore $\F\cap\Sigma$ is a s.r.f.s. $\fimdem$

%%%%%%%%% OS PONTOS SINGULARES NA GEODESICA OU SAO ISOLADOS OU SAO TODA A GEODESICA
%%%%%%%%%%
\begin{prop}
\label{pontos-singulares-na-geodesica}
Let $\F$ be a s.r.f.s. on a complete riemannian manifold $M$ and $\gamma$ a geodesic orthogonal to the leaves of   $\F$. Then the singular points of $\gamma$ are either all the points of  $\gamma $ or isolated points of $\gamma.$ 
\end{prop}
\dem  Since the set of regular points on $\gamma$ is open, we can suppose that $q=\gamma(0)$ is a singular point and that $\gamma(t)$ is a regular point for $-\delta<t<0.$ We shall show that there exists $\epsilon>0$ such that $\gamma(t)$ is also a regular point for $0<t<\epsilon.$

At first we note that we can choose $t_{0}<0$ such that $q$ is a focal point of $L_{\gamma(t_{0})}.$ To see this let $\tub(P_{q})$ be a tubular neighborhood of a plaque $P_{q}$ and $t_{0}<0$ such that $\gamma(t_{0})\in \tub(P_{q}).$ Since $L_{\gamma(t_{0})}$ is a regular leaf and $q$ is a singular point, it follows from item c) of  Proposition \ref{slice-uniao-secoes} that $L_{\gamma(t_{0})}\cap \Sigma_{q}$ is not empty. Then we can join this submanifold to $q$ with  geodesics that belong to $\Sigma_{q}.$ Since $\F$ is a riemannian foliation, these geodesics are also orthogonal to $L_{\gamma(t_{0})}\cap \Sigma_{q}.$ This implies that $q$ is a focal point. 

 Since focal points are isolated along $\gamma,$ we can choose $\epsilon>0,$ such that $\gamma(t)$ is not a focal point of $P_{\gamma(t_{0})}$ along $\gamma$ for $0<t<\epsilon.$ 

Suppose there exists $0<t_{1}<\epsilon$ such that $x=\gamma(t_{1})$ is a singular point. Let $\sigma$ a local section that contains $\gamma(t_{0}).$ 
Let $U$ an open set of $ \nu_{x}L$ such that $\widetilde{\Sigma}_{x}:=\exp_{x}(U)$ contains $\gamma(t_{0})$ and is contained in a convex neighborhood of $x.$ We note that $\widetilde{\Sigma}_{x}$ is not contained in a tubular neighborhood of $P_{x}$ and hence is not a slice. 

We have:
\[\sigma\subset\widetilde{\Sigma}_{x}.\]
Since $x$ is a singular point, we have:
\[\dim \sigma<\dim \widetilde{\Sigma}_{x}.\]

The equations above implies that $\dim  P_{\gamma(t_{0})}\cap\widetilde{\Sigma}_{x}>0.$ Hence we can find geodesics in $\widetilde{\Sigma}_{x}$ that join $x$ to the submanifold  $ P_{\gamma(t_{0})}\cap\widetilde{\Sigma}_{x}.$ Since the foliation is a riemannian foliation, these geodesics are also orthogonal to $ P_{\gamma(t_{0})}\cap\widetilde{\Sigma}_{x}$ and hence $x$ is a focal point of this submanifold. This contradicts our choice of $\epsilon$ and completes the proof. $\fimdem$

%%%%%%%%%%%%%%%%%%%%%%%%%%%%% LEMA HOL%%%%%%% 

In what follows we shall  need  a result of Heintze, Liu and Olmos.

\begin{prop}[Heintze, Liu and Olmos \cite{HOL}]
\label{lema-HOL}
Let $M$ be a complete riemannian manifold, $L$ be an immersed submanifold of  $M$ with  globally flat normal bundle and  $\xi$ be a normal parallel field along $L.$ Suppose that  $\sigma_{x}:=\exp_{x}(\nu_{x}L)$ is a totally geodesic complete submanifold for all $x\in L,$ that means, $L$ has sections. Then 
\begin{enumerate}
\item $ d \eta_{\xi}(v)$ is orthogonal to $\sigma_{x}$ at $\eta_{\xi}(x)$ for all $v\in T_{x}L.$
\item Suppose that $p$ is not a critical point of the map $\eta_{\xi}.$ Then there exists  a neighborhood $U$ of $p$ in $L$ such that  $\eta_{\xi}(U)$ is an embedded submanifold, which meets $\sigma_{x}$ orthogonally
and has globally flat normal bundle.
 In addition,  a parallel normal field along $U$ transported to $\eta_{\xi}(U)$ by parallel translation along the geodesics $\exp(t \xi)$ is a parallel normal field along $\eta_{\xi}(U).$
\end{enumerate}
\end{prop} 
 Let $\exp^{\perp}$  denote the restriction of $\exp$ to $\nu(L).$ We recall that for each $w\in T_{\xi_{0}}\nu(L)$ there exists only one $w_{\mathsf{t}}\in T_{\xi_{0}}\nu(L)$ ( the tangential vector) and one $w_{\mathsf{n}} \in T_{\xi_{0}}\nu(L)$ ( the normal vector) such that
\begin{enumerate}
\item $w=w_{\mathsf{t}}+w_{\mathsf{n}},$
\item $d\Pi(w_{\mathsf{n}})=0,$ where $\Pi: \nu(L)\rightarrow L$ is the natural projection,
\item$ w_{\mathsf{t}}=\xi^{'}(0),$  where $\xi(t)$ is the normal parallel field with $\xi(0)=\xi_{0}.$
\end{enumerate}
We also recall that $z=\exp^{\perp}(\xi_{0})$ is a focal point with multiplicity $k$ along $\exp^{\perp}(t\,\xi_{0})$ if and only if $\dim \ker d \exp^{\perp}_{\xi_{0}}=k.$ We call $z$ a focal point of $L$  of \emph{tangencial type} if $\ker d \exp^{\perp}_{\xi_{0}}$ only consists of tangential vectors. 

\begin{cor}
\label{ponto-focal-proximo}
Let $L$ be a submanifold defined as above, $p\in L$ and $\xi_{0}\in\nu_{p}L.$ Suppose that  the  point $z=\exp_{p}(\xi_{0})$ is a focal point of $L$ along $\exp_{p}(t\,\xi_{0})$ that belongs to a normal neighborhood of $p.$ Then $z$ is a focal point of tangencial type.
\end{cor}
\dem  If $z$ is a focal point, then there exists $w\in T_{\xi_{0}}\nu(L)$ such that $\|d \exp_{\xi_{0}}^{\perp}(w)\|$ $=0.$ It follows from the above proposition  that
\[< d\exp^{\perp}_{\xi_{0}}w_{\mathsf{n}},d \exp^{\perp}_{\xi_{0}}w_{\mathsf{t}}>_{\exp^{\perp}(\xi_{0})}=0\]
and hence $\|d \exp^{\perp}_{\xi_{0}}(w_{\mathsf{n}})\|=0.$ Since $z$ belongs to a normal neighborhood,  $w_{\mathsf{n}}$ must be zero.  We conclude that $w=w_{\mathsf{t}}.$ $\fimdem$

Now we can show one of our main results.

%%%%%%%%%%%%%%%%%%%%%FOLHAS DE FRSS SAO LOCALMENTE EQUIFOCAIS%%%%%%

\begin{teo}
\label{frss-eh-equifocal}
Let $\F$  be a singular riemannian foliation with sections on a complete riemannian manifold  $M.$ Then the regular leaves are locally equifocal. In addition, if  all the leaves are compact, then the union of  regular leaves that are equifocal is an open and dense set in $M.$
\end{teo}

To prove it, we need  the following lemma.

\begin{lema}
\label{lema-teo-frss-sao-localmente-equifocais}
Let $\tub(P_{q})$ be a tubular neighborhood of a plaque $P_{q},$
$x_{0}\in \tub(P_{q}),$   and  $\xi\in \nu P_{x_{0}}$ such that  $\exp_{x_{0}}(\xi)=q.$  We also suppose that $q$ is the only singular point on the segment of geodesic $\exp_{x_{0}}(t\,\xi)\cap \tub(P_{q}).$  Then we can find a neighborhood $U$ of $x_{0}$ in $P_{x_{0}}$ with the following properties:
\begin{enumerate}
\item[1)] $\nu U$ is globally flat  and we can define the parallel normal field $\xi$ on $U.$  
\item[2)] There exists a number $\epsilon>0$ such that,  
for each $x\in U,$  $\gamma_{x}\subset \tub(P_{q}),$ where  $\gamma_{x}(t):=\exp_{x}(t\,\xi)$ and $t\in [-\epsilon,1+\epsilon].$ 
\item[3)] The regular points of the foliation  $\F|_{\tub(P_{q})}$ are not critical values of the maps $\eta_{t\,\xi}|_{U}.$
\item[4)] $\eta_{t\,\xi}(U)\subset L_{\gamma_{x_{0}}(t)}.$ 
\item[5)] $\eta_{t\,\xi}:U\rightarrow \eta_{t\,\xi}(U)$ is a  local diffeomorphism  for $t\neq 1.$
\item[6)] $\dim \rank D \eta_{\xi}$ is constant on $U.$ 
\end{enumerate}
\end{lema}

\dem The item 1) follows from the fact that $\F$ has sections and one can show 2) with standard arguments. 

3) Let $p=\eta_{r\,\xi}(x_{1})$ be a regular point of the foliation and suppose that $x_{1}$ is a critical point of the map $\eta_{r\,\xi}|_{U}.$ Then there exists a Jacobi field $J(t)$ along the geodesic $\gamma_{x_{1}}$ such that $J(r)=0.$ In particular there exists a smooth curve $\beta(t)\subset P_{x_{0}}$ such that  $J(t)=\frac{\partial}{\partial s}\exp_{\beta(s)}(t\,\xi)$ and $\beta(0)=x_{1}.$  

Since focal points are isolated along $\gamma_{x_{1}}(t),$ there exists a regular point of the foliation $\tilde{p}=\gamma_{x_{1}}(\tilde{r})$ that is not a focal point of $P_{x_{1}}$ along $\gamma_{x_{1}}.$  It follows from  Proposition \ref{lema-HOL} that there exists a neighborhood $V$ of $x_{1}$ in $P_{x_{0}}$ such that the embedded submanifold  $\eta_{\tilde{r}\,\xi}(V)$ is orthogonal to the sections that it meets. 
Hence $\eta_{\tilde{r}\,\xi}(V)$ is tangent to the plaques near to $ P_{\tilde{p}}.$ Since $\eta_{\tilde{r}\,\xi}(V)$ has the dimension of the regular leaves, $\eta_{\tilde{r}\,\xi}(V)$ is an open subset of  $P_{\tilde{p}}.$

Since  we can choose $\tilde{p}$ so close to $p$ as necessarily, we can suppose that $p$ and $\tilde{p}$ belong to a neighborhood $W$ that contains only regular points of the foliation and such that $\F|_{W}$ are pre image of an integrable riemannian submersion $\pi:W\rightarrow B.$  It follows from  Proposition \ref{lema-HOL} that $\gamma_{\beta(s)}^{'}(\tilde{r})$ is a  parallel field along the curve $\eta_{\tilde{r}\xi}\circ\beta(s) \subset \eta_{\tilde{r}\,\xi}(V)\subset P_{\tilde{p}}.$ Therefore $\gamma_{\beta(s)}(t)\cap W$ are horizontal lift of a geodesic in $B$ ( the basis of the riemannian submersion $\pi$).  This implies that $J(r)\neq 0$ This contradicts the assumption that $p$ is a focal point and completes the proof of item 3).

4) At first we check the item 4) for each $t\neq 1.$ Fix a $t_{0}\neq 1$ and define $K:=\{k\in U$ such that $\eta_{t_{0}\,\xi}(k)\in P_{\gamma_{x_{0}}(t_{0})}\}.$ Since $\gamma_{x_{0}}(t_{0})$ is a regular point of the foliation, it follows from the item 3) that all the points of $P_{\gamma_{x_{0}}(t_{0})}$ are regular values of the map $\eta_{t_{0}\,\xi}.$ Hence for each $k\in K$ there exists a neighborhood $V$ of $k$ in $U$ such that $\eta_{t_{0}\,\xi}(V)$ is an embedded submanifold. As we have note in the proof of item 3), $\eta_{t_{0}\,\xi}(V)$ is an open set of $P_{\gamma_{x_{0}}(t_{0})},$ because this embedded submanifold is orthogonal to the sections and has the same dimension of the plaques. We conclude then that $K$ is an open set. One can prove that  $K$ is closed using standard arguments and the fact that the plaques are equidistant. Since $U$ is connected, $K=U.$

Now we check the item 4) for $t=1.$ We define  $f(x,t):=d(\eta_{t\,\xi}(x),P_{q})- d(\eta_{t\,\xi}(x_{0}),P_{q}).$ As we have seen above   $\eta_{t\,\xi}(x)$ and $\eta_{t\,\xi}(x_{0})$  belong to the same plaque, for $t\neq 1.$ This means that $f(x,t)=0$ for all  $t\neq 1$ and hence  $f(x,1)=0,$ i.e.,  $\eta_{\xi}(x)\subset P_{q}.$

5) The item 5) follows from the item 3) and 4).

6) Fix a point $x_{1}\in U.$ It follows from  Corollary  \ref{ponto-focal-proximo} that the focal points of $U$ along $\gamma_{x}(t)$ are of tangential type. This means that  $\gamma_{x}(t_{0})$  is a focal point of $U$ along $\gamma_{x}$ with multiplicity $k$ if and only if $x$ is a critical point of $\eta_{t_{0}\,\xi}$ and  $\dim\ker d\eta_{t_{0}\,\xi}(x)=k.$   
In addition, it follows from the item 5) that 
the map $\eta_{t\,\xi}$ might not be a diffeomorphism only for
$t=1.$   Therefore we have
\begin{eqnarray}
\label{pt-focal-gamma1} 
m (\gamma_{x})= \dim\ker d\eta_{\xi}(x),
\end{eqnarray}
 where $m(\gamma_{x})$ denote the number of focal points on $\gamma_{x}(t),$ each counted with its multiplicities.

On the other hand, we have  
\begin{eqnarray}
\label{morse-desigualdade-ptfocais}
m (\gamma_{x})\geq m(\gamma_{x_{1}})
\end{eqnarray}
for $x$ in a neighborhood of $x_{1}$ in $U.$ Indeed one can argue like 
Q.M. Wang\cite{Wang1} to see that equation \ref{morse-desigualdade-ptfocais} follows from the Morse index theorem. 

Equations  (\ref{pt-focal-gamma1}) and (\ref{morse-desigualdade-ptfocais})  together with the elementary expression \linebreak  $\dim\ker d\eta_{\xi}(x)\leq\dim  \ker d\eta_{\xi}(x_{1})$ imply that  $\dim \ker d\eta_{\xi}$  is constant on a neighborhood of $x_{1}$ in $U.$ Since this hold for each $x_{1}\in U,$ we conclude that that  $\dim \ker d\eta_{\xi}$ is constant on $U.$ $\fimdem$

 \textbf{Proof of  Theorem \ref{frss-eh-equifocal}} 

Let $L$ be a leaf of $\F,$ $U$ be a open set of $L$ that has normal bundle globally flat and $\xi$ a parallel normal field along $U.$ At first, we shall prove that $\dim\rank d\eta_{\xi}|_{U}$ is constant, i.e., $L$ is locally equifocal. 

Let $p\in U.$ Since singular points are isolated along 
   $\gamma_{p}(t)=\exp_{p}(t\,\xi)|_{[-\epsilon,1+\epsilon]}$ (see Proposition \ref{pontos-singulares-na-geodesica}), 
we can cover this arc of geodesic with a finite number of tubular neighborhood  $\tub(P_{\gamma_{p}(t_{i})})$ , where $t_{0}=0$ and  $t_{n}=1.$   

Let  $P_{\gamma_{p}(r_{i})}$ be regular plaques that belong to  $\tub(P_{\gamma_{p}(t_{i-1}}))\cap \tub(P_{\gamma_{p}(t_{i})}),$ where $t_{i-1}<r_{i}<t_{i}.$ Applying the Lemma  \ref{lema-teo-frss-sao-localmente-equifocais} and  Proposition  \ref{lema-HOL}, we can find  an open set $U_{0}\subset P_{p},$ of the plaque $P_{p},$    open sets $U_{i}\subset P_{\gamma_{p}(r_{i})}$ of the plaques $P_{\gamma_{p}(r_{i})}$ and  parallel normal fields  $\xi_{i}$ along  $U_{i},$ with the following properties: 
\begin{enumerate}
\item[1)] For each  $U_{i},$ the parallel normal field  $\xi_{i}$  is tangent to the geodesics $\gamma_{x}(t),$ where $x\in U_{0};$
\item[2)] $ \eta_{\xi_{i}}:U_{i}\rightarrow U_{i+1}$ is a  local diffeomorphism for $i<n;$
\item[3)] $\eta_{\xi}|_{U_{0}}=\eta_{\xi_{n}}\circ\eta_{\xi_{n-1}}\circ\cdots\circ\eta_{\xi_{0}}|_{U_{0}}$ 
\end{enumerate}
   
Since $\dim\rank d\eta_{\xi_{i}}$ is constant on $U_{i}$,  $\dim d\eta_{\xi}$ is constant on $U_{0}.$ Since this hold for each $p\in U,$ 
$\dim d\eta_{\xi}$ is constant on $U,$ i.e., $L$ is locally equifocal. 

 At last we check what happens when the leaves of the foliation are compact. According to Molino  (Proposition 3.7, page 95 \cite{Molino}), the union of the regular leaves with trivial holonomie of a singular riemannian foliation is an open and dense set in the set of the regular points. In addition, the set of regular points is an open and dense set  in $M$ (see page 197 \cite{Molino}). Since the leaves of a s.r.f.s. that have trivial holonomie are exactly the leaves that have normal bundle globally flat,  the union of  regular leaves that are equifocal is an open and dense set in $M.$ $\fimdem$
  
%%%%%%%%%%%%%%%%%%COROLARIO: FOLHAS PARALELAS PERTENCEM A FOLHEACAO
\pagebreak

\begin{cor}
\label{cor-map-paralelo}
Let $\F$ be a s.r.f.s. on a complete riemannian manifold  $M$ and $L$ be a regular leaf of  $\F.$ 
\begin{enumerate}
\item[a)] Let  $\beta(s)\subset L$ be a smooth curve  and  $\xi$ a parallel normal field along  $\beta(s).$ Then the curve $\exp_{\beta(s)}(\xi)$ belongs to a leaf of the foliation.  
\item[b)] Let  $L$ be a regular equifocal leaf and $\Xi$ denote the set of all parallel normal fields along $L.$ Then $\F=\{\eta_{\xi}(L)\}_{\xi\in \, \Xi}.$ 
\end{enumerate}
\end{cor}

\dem a) The item a) can be easily proved using the item 4) of Lemma \ref{lema-teo-frss-sao-localmente-equifocais} and gluing tubular neighborhoods as we have already done in the proof of Theorem 
\ref{frss-eh-equifocal}.

b) \afir Let $x_{0}\in L$ and $q=\eta_{\xi}(x_{0}).$ Then there exists a neighborhood $U\subset L$ of $x_{0}$ in $L$ such that $\eta_{\xi}(U)\subset P_{q}$ is an open set in $L_{q}.$

To check this statement is enough to suppose that $x_{0}\in\tub(P_{q}),$ for the general case can be proved gluing tubular neighborhoods as we have done in proof of  Theorem  \ref{frss-eh-equifocal}. Now the statement follows if we note that $\eta_{\xi}:U\rightarrow P_{q}$ is a submersion whose fibers are the intersections of $U$ with the slices of $P_{q}.$ To see that each fiber is contained in  $U\cap \Sigma$ one can use the fact that the rank of $d\eta_{\xi}$ is constant  and the fact that the foliation is riemannian. To see that each fiber contains $U\cap \Sigma$ one can use  item 4) of Lemma \ref{lema-teo-frss-sao-localmente-equifocais} together with item a) of   Proposition \ref{slice-uniao-secoes}.

It follows from the above statement that $\eta_{\xi}(L)$ is an open set in $L_{q}.$ We shall see that $\eta_{\xi}(L)$ is also a closed set in $L_{q}.$ Let $z\in L_{q}$ and $\{z_{i}\}$ a sequence in $\eta_{\xi}(L)$ such that $z_{i}\rightarrow z.$

At first suppose  that $L_{q}$ is a regular leaf. If follows from  Proposition \ref{lema-HOL} that there exists a parallel normal field $\hat{\xi}$ along $\eta_{\xi}(L)$ such that $\hat{\xi}_{\eta_{\xi}(x)}$ is tangent to the geodesic $\exp_{x}(t\,\xi).$ Since the normal bundle of $P_{z}$ is globally flat, we can extend $\hat{\xi}$ along $P_{z}.$ The item 4) of  Lemma \ref{lema-teo-frss-sao-localmente-equifocais} implies that $\eta_{-\hat{\xi}}:P_{z}\rightarrow L.$ By construction $\eta_{\xi}\circ\eta_{-\hat{\xi}}(z_{i})=z_{i}.$ Therefore $\eta_{\xi}\circ\eta_{-\hat{\xi}}(z)=z.$  This means that $z\in \eta_{\xi}(L).$

At last suppose that $L_{q}$ is a singular leaf. There exists $x_{i}\in L$ such that $z_{i}=\eta_{\xi}(x_{i})\in P_{z}.$ We can find a $s<1$ such that $y_{i}=\eta_{s\,\xi}(x_{i})$ is a regular point. Since $y_{i}$ is a regular point, the plaque $P_{y_{i}}$ is an open set of $\eta_{s\,\xi}(L)$ as we have proved above. There exists a parallel normal field $\hat{\xi}$ along $P_{y_{i}}$ such that $\eta_{\hat{\xi}}\circ\eta_{s\,\xi}=\eta_{\xi}.$ 
 
It follows from item 4) Lemma \ref{lema-teo-frss-sao-localmente-equifocais} that $\eta_{\hat{\xi}}(P_{y_{i}})\subset P_{z}.$ On the other hand, since the foliation is singular,  the plaque $P_{y_{i}}$ intercept the slice $\Sigma_{z}.$ These two facts  imply that $z\in  \eta_{\hat{\xi}}(P_{y_{i}}).$ Therefore $z\in \eta_{\xi}(L).$ $\fimdem$

%%%%%%%%%%%%%%TEOREMA DO SLICE%%%%%%%%%%%%%%%%%%%%%%
%%%%%%%%%%%%%%%%%%%%%%%%%%%%%%%%%%%%%%%%%%%%%%%%%%%%%%

Let $\F$ be a foliation on a manifold $M^{n},$ $\widetilde{\F}$  a foliation on a manifold $\widetilde{M}^{n}$ and $\varphi: M\rightarrow \widetilde{M}$ a diffeomorphism. We say that $\varphi$ is a \emph{diffeomorphism between $\F$ and $\widetilde{\F}$} if each leaf $L$ of $\F$ is diffeomorphic to a leaf $\widetilde{L}$ of $\widetilde{\F}$.

\begin{teo}[slice theorem]
\label{sliceteorema}
Let $\F$ be a singular riemannian foliation with sections on a complete riemannian manifold $M$  and  $\Sigma_{q}$ a slice at a point  $q\in M.$ Then  $\F$ restrict to  $\Sigma_{q}$ is diffeomorphic to an isoparametric foliation on an open set of $\mathbf{R}^{n},$  where $n$ is the dimension of   $\Sigma_{q}.$ 
\end{teo}

\dem According to H. Boualem (see Proposition 1.2.3 and Lemma 1.2.4 of \cite{Boualem}), we have
\begin{enumerate}
\item[1)] the map $\exp^{-1}$ is a diffeomorphim between   the foliation $\F|_{\Sigma_{q}}$ and a singular riemannian foliation with sections  $\widetilde{F}$ on an open set of the inner product space  $(T_{q}\Sigma_{q},<,\,>_{q}),$ where  $<,\,>_{q}$ denote the metric of $T_{q}M,$ 
\item[2)] the sections of the singular  foliation $\widetilde{F}$ are the vector subspaces $\exp^{-1}(\sigma),$ where $\sigma$ are the local section of $\F.$ 
\end{enumerate}

Let $<,>_{0}$ denote the canonical euclidian product. Then there exists a positive  definite symmetric matrix $A$ such that $<X,Y>_{q}=<A\,X,Y>_{0}$  The isometry $\sqrt{A}:(T_{q}\Sigma_{q},<,>_{q})\rightarrow (\mathbf{R}^{n},<,>_{0})$ is a diffeomorphismus between  the foliation $\widetilde{F}$ and a singular riemannian foliation with section $\widehat{F}$ on an open set of the inner product space $(\mathbf{R}^{n},<,>_{0}).$ Since $0\in \mathbf{R}^{n}$ is a singular leaf of the foliation $\widehat{\F},$ the leaves of this foliation belong to spheres in the euclidian space.

\afir 1: The restriction the foliation $\widehat{\F}$ to a sphere $\mathbf{S}^{n-1}(r)$  is a singular riemannian foliation with sections on $\mathbf{S}^{n-1}(r).$

The first step to check this statement is to note that $\widehat{\F}|_{\mathbf{S}^{n-1}(r)}$  is a singular foliation, for $\widehat{\F}$ is a singular foliation. Next we have to note that if $\widehat{\sigma}$ is a section of $\widehat{F}$ then $\sigma_{s}:=\widehat{\sigma}\cap \mathbf{S}^{n-1}(r)$ is a section of the foliation  $\widehat{\F}|_{\mathbf{S}^{n-1}(r)}.$ To conclude, we have to note that $\widehat{\F}|_{\mathbf{S}^{n-1}(r)}$ is a transnormal system.  Let $\gamma$ be a geodesic of $\mathbf{S}^{n-1}(r)$  that
is orthogonal to a leaf $L_{\gamma(0)}$ of $\widehat{\F}|_{\mathbf{S}^{n-1}(r)}.$ Since a slice of $\widehat{F}$ is a union of sections, $\gamma$ is tangent to a section $\hat{\sigma}$ at the point $\gamma(0)$ and hence is tangent to a section $\sigma_{s}$ of $\widehat{\F}|_{\mathbf{S}^{n-1}(r)}$ at the point $\gamma(0).$  
This implies that $\gamma\subset\sigma_{s},$ which means that $\gamma$ is orthogonal to each leaf that it meets, i.e.,  the partition is transnormal. 

Now  Theorem \ref{frss-eh-equifocal} guarantees that the leaves of a singular riemannian foliation with sections are locally equifocal. Therefore the leaves of $\widehat{\F}|_{\mathbf{S}^{n-1}(r)}$ are locally equifocal. 

The next statement  follows from standard calculations on space forms. 

\afir 2: The locally equifocal submanifolds in $\mathbf{S}^{n-1}(r)$ are isoparametric submanifold in  $\mathbf{S}^{n-1}(r).$

Since isoparametric submanifold in spheres are isoparametric submanifolds in euclidian spaces (see Palais and Terng, Proposition 6.3.17 \cite{PTlivro}), we can conclude that the regular leaves of $\widehat{\F}$ are isoparametric submanifold in an open set of the euclidian space $\mathbf{R}^{n}.$

At last, we note that Corollary  \ref{cor-map-paralelo}
 implies that the singular leaves of $\widehat{\F}$ are the focal leaves. Therefore $\widehat{\F}$ is an isoparametric foliation on an open set of the euclidian space and this completes the proof of the theorem. $\fimdem$

%%%%%%%%%%%%%%%%%COR: ESTRATIFICACAO SINGULAR
\begin{cor}
\label{estratificacao-singular}
Let $\F$ be a  s.r.f.s. on a complete riemannian manifold $M$ and $\sigma$ be a local section contained in a slice  $\Sigma_{q}$ of dimension $n.$ According to the slice theorem there exist an open set $U\subset \mathbf{R}^{n}$ and a diffeomorphism  $\Psi:\Sigma_{q}\rightarrow U$ that sends  the foliation $\F\cap \Sigma_{q}$  to  an isoparametric foliation $\widehat{\F}$ on $U.$ 
Then the set of singular points of $\F$ contained in $\sigma$ is a finite union of totally geodesic hypersurfaces that are sent by $\Psi$ onto the focal hyperplanes of $\widehat{\F}$ contained in a section of this isoparametric foliation.   
\end{cor}
We shall call \emph{singular stratification of the local section} $\sigma$ this set of singular points of $\F$ contained in $\sigma.$  

\dem It follows from Molino \cite{Molino}(page 194, Proposition 6.3) that the intersection of the singular leaves with a section is a union of totally geodesic submanifolds. Now the slice theorem implies that these totally geodesic submanifolds are in fact hypersurfaces  that are diffeomorphic to focal hyperplanes. $\fimdem$

%%%%%%%%%%%%%%%%%%%%%%%%%%%%%%%%%SISTEMA LOCAL DE COORDENADAS%%%%%%%%%%%%%%
%%%%%%%%%%%%%%%%%%%%%%%%%%%%%%%%%%%%%%%%%%%%%%%%%%%%%%%%%%%%%%%%%%%%%%%%%%%%%%%%
\begin{prop}
\label{frss-transnormal}
Let $\F$ be a s.r.f.s. on a complete riemannian manifold  $M$ and $q\in M.$ Then there exist a tubular neighborhood  $\tub(P_{q}),$ an open set $W\subset \mathbf{R}^{k}$  and a transnormal map   $H : \tub(P_{q})\rightarrow W$ such that the preimages  of  $H$ are leaves  of the singular foliation  $\F|_{\tub(P_{q})}.$ The leaf  $H^{-1}(c)$ is regular if and only if $c$ is  a regular value. 
\end{prop}
\dem We start recalling a result that can be found in the book of Palais and Terng.

\begin{lema}[Theorem 6.4.4. page 129 \cite{PTlivro}]
Let $N$ be a rank $k$ isoparametric submanifold in $\textbf{R}^{n},$ $W$ the associated Coxeter group, $q$ a point on $N,$  $ \nu_{q}=q+\nu(N)$ the affine normal plane at $q$ and $u_{1},\cdots, u_{k}$ be a set of generators of the $W$-invariant polynomials on $\nu_{q}.$ Then $u=(u_{1},\cdots, u_{k})$ extends uniquely to an isoparametric polynomial map $g:\mathbf{R}^{n}\rightarrow \textbf{R}^{k}$ having $N$ as a regular level set. Moreover, 
\begin{enumerate} 
\item[1)] each regular set is connected,
\item[2)] the focal set of $N$ is the set of critical points of $g,$
\item[3)] $\nu_{q}\cap N=W\cdot q,$
\item[4)] $g(\textbf{R}^{n})=u(\nu_{q}),$
\item[5)] for $x\in\nu_{q},$ $g(x)$ is a regular value if and only if $x$ is $W$-regular,
\item[6)] $\nu(N)$ is globally flat.
\end{enumerate}
\end{lema}

The above result implies that the leaves of the isoparametric foliation, which  has $N$ as a leaf, can be described as pre image of a map $g.$ Note that this is even true if $N$ is not a full isoparametric submanifold of $\textbf{R}^{n}.$

Now we define $\widetilde{H}:\Sigma_{q}\rightarrow \textbf{R}^{k}$ as $\widetilde{H}:=g\circ\Psi,$ where $\Psi:\sigma_{q}\rightarrow \textbf{R}^{n}$ is the diffeomorphism given by the slice theorem that sends $\F|_{\Sigma_{q}}$ to an isoparametric foliation on an open set $W$ of $\textbf{R}^{n}.$

Since $\F$ is a singular foliation, there exists a projection $\Pi:\tub(P_{q})\rightarrow\Sigma_{q}$ such that $\Pi(P)=P\cap\Sigma_{q}$ for each plaque $P.$

Finally we define $H:=\widetilde{H}\circ\Pi.$ Then the preimages of $H$ are leaves of the foliation $\F|_{\tub(P_{q})}.$ 

The statement below, which can be found in Molino \cite{Molino}[page 77], implies that $H$ is a transnormal map.

\afir Let $U$ a simple neighborhood of a riemannian foliation ( with section) and $H:U\rightarrow \widetilde{U}\subset\textbf{R}^{k}$ such that $H^{-1}(c)$ are leaves of $\F|_{U}.$ Then we can choose a metric for $\widetilde{U}$ such that $H:U\rightarrow \widetilde{U}$ is a (integrable)  riemannian submersion. $\fimdem$

%%%%%%%%%%%%%%%%%
%%%%%%%%%%%%%%%
%%%%%%%%%%%%%%%

% second part of the article: singular holonomy and new examples
\section{Singular Holonomy} 

The slice theorem give us a description of  the plaques of a singular riemannian foliation with sections. However, it doesn't assure us if two different plaques belong to the same leaf.
To get such kind of information, we must extend the concept of holonomy to describe not only what happens near a regular leaf but also what happens in a neighborhood of a singular leaf. 

In this section, we shall introduce the concept of singular holonomy and establish some of its properties.   

\begin{prop}
\label{prop-holonomia-singular}
Let $\F$ be a s.r.f.s. on a complete riemannian manifold $M,$ $L_{p}$ a regular leaf, $\sigma$ a local section and $\beta(s)\subset L_{p}$ a smooth curve, where $p=\beta(0)$ and $\beta(1)$ belong to $\sigma.$ Let $[\beta]$ denote the homotopy class of $\beta.$ Then there exists an isometry $\varphi_{[\beta]}:U\rightarrow W,$ where the source $U$ and target $W$ contain $\sigma,$ which has the following properties:
\begin{enumerate}
\item[1)]$\varphi_{[\beta]}(x)\in L_{x}$ for each $x\in\sigma,$
\item[2)]$d\varphi_{[\beta]}\xi(0)=\xi(1),$ where $\xi(s)$ is a  parallel normal field along $\beta(s).$
\end{enumerate}
\end{prop}
\dem Since $\sigma$ is a local section, for each $x\in\sigma,$ there exists only
one  $\xi\in T_{p}\sigma$ such that $\exp_{p}(\xi)=x.$ Let $\xi(t)$ be the parallel transport of $\xi$ along $\beta$ and define $\varphi_{\beta}(x):=\exp_{\beta(1)}(\xi(1)).$ It's easy to see that $\varphi_{\beta}$ is a bijection. It follows from  Corollary \ref{cor-map-paralelo} that $\exp_{\beta}(\xi)\subset F_{x}$ and this proves a part of item 1. Since $\varphi_{\beta}$ is an extension of the holonomy map, $ d\varphi_{\beta}\xi(0)=\xi(1),$ and this proves a part of item 2.  The fact that $\varphi_{\beta}$ is an extension of the holonomy map implies that the restriction of $\varphi_{\beta}$ to a small neighborhood of $\sigma$ depend only on the homotopy class of $\beta.$ Since isometries are determined by the image of a point and the derivative at this point, is enough to prove that $\varphi_{\beta}$ is an isometry to see that $\varphi_{\beta}$ depends only of the homotopy class of $\beta.$
To see that $\varphi_{\beta}$ is an isometry it's enough to check the following statement.

\afir Given a point $x_{0}\in\sigma$ there exists an open set $V\subset\sigma$ of $x_{0}$ in $\sigma$  such that $d(x_{1},x_{0})=d(\varphi_{\beta}(x_{1}),\varphi_{\beta}(x_{0})),$ for each $x_{1}\in V.$

To check the statement let $\xi_{0}(s)$ and $\xi_{1}(s)$ be normal parallel fields along $\beta(s)$ such that $x_{j}=\exp_{p}(\xi_{j}(0))$ for $j=0,1.$ Define $\alpha_{j}(s)=\exp_{\beta(s)}(\xi_{j}(s))$ for $j=0,1.$ Since $\varphi_{\beta}(x_{j})=\alpha_{j}(1)$ the statement follows from the following equation
\[d(\alpha_{0}(s),\alpha_{1}(s))=d(\alpha_{0}(0),\alpha_{1}(0))\]
and this equation follows from the following facts:
\begin{enumerate} 
\item $\alpha_{j}(s)\in L_{x_{j}}$
\item singular riemannian foliations are locally equidistant,
\item $\alpha_{0}(s)$ and $\alpha_{1}(s)$ are always in the same local section. $\fimdem$
\end{enumerate}

\begin{definicao}
\emph{
The pseudosubgroup of isometries  generated by the isometries constructed above is called \emph{pseudogroup  of singular holonomy of the local section $\sigma.$} Let $\Holsing(\sigma)$ denote this pseudogroup.
}
\end{definicao}
\begin{prop}
Let $\F$ be a s.r.f.s. on a complete riemannian manifold $M$ and $\sigma$ a local section. Then the reflections in the hypersurfaces of the singular stratification of the local section $\sigma$ let $\F\cap\sigma$  invariant. Moreover these reflections are elements of $\Holsing  (\sigma).$
\end{prop}
\dem The proposition is already true if the singular foliation is an isoparametric foliation on an euclidean space. In what follows we shall use this fact and the slice theorem to construct the desired reflections.

Let $S$ be a complete totally geodesic hypersurface of the singular stratification of the local section $\sigma$ and $\Sigma$ be a slice of a point of $S$ and hence that contains $\sigma.$  It follows from the slice theorem that there exists a diffeomorphism $\Psi:\Sigma\rightarrow V\subset\mathbf{R}^{n}$ that sends $\F\cap\Sigma$ to an isoparametric foliation $\widetilde{F}$ on an open set $V$ of $\mathbf{R}^{n}.$ Let $p\in \sigma$ be a regular point, $\tilde{L}:=\Psi(L_{p}\cap\Sigma)$ and $\tilde{\sigma}:=\Psi(\sigma).$ We note that $\tilde{\sigma}$ is a local section of the isoparametric foliation $\widetilde{F}.$ 

It follows from  Corollary \ref{estratificacao-singular} and from the theory of isoparametric submanifolds \cite{PTlivro} that $\tilde{S}:=\Psi(S)$ is a focal hyperplane associated to a curvature distribution $E.$
Let $\beta\subset \Sigma\cap\F$ with $\beta(0)=p$ and $\beta(1)\in\sigma$  such that $\tilde{\beta}:=\Psi\circ\beta$ is tangent to the distribution $E.$
Finally let $z\in S,$ $\xi\in T_{p}\sigma$  such that $\exp_{p}(\xi)=z$ and $\xi(s)$ the parallel transport of $\xi$ along $\beta.$ 

\afir $\exp_{\beta(s)}(\xi)=z.$

To check this statement, we recall  that $\tilde{S}\subset\tilde{\sigma}_{\tilde{\beta}},$ where $\tilde{\sigma}_{\tilde{\beta}}$ i a local sectionn of $\tilde{\F}$ that contains $\tilde{\beta(s)}$(see Theorem 6.2.9 \cite{PTlivro}). Therefore $S\subset\sigma_{\beta(s)}.$ On the other hand, it follows from  Corollary \ref{cor-map-paralelo}
 that $\exp_{\beta(s)}(\xi)\subset P_{z}.$  Hence $\exp_{\beta(s)}(\xi)\subset P_{z}\cap S.$ Now the statement follows from the fact that $P_{z}\cap S=\{z\}.$ 

This statement implies that the isometry $\varphi_{[\beta]}$ let the points of $S$ fixed. Therefore $\varphi_{[\beta]}$ is a reflection in a totally geodesic hypersurface. Since $\varphi_{[\beta]}(x)\in L_{x},$  these reflections let $\F\cap\sigma$  invariant. $\fimdem$  

\pagebreak

\begin{prop}
Let $\F$ be a s.r.f.s. on a complete riemannian manifold $M.$ Suppose that the leaves are compact and that the holonomies of regular leaves are trivial. Let $\sigma$ be a local section and $\Omega$ a connected component of the set obtained removing the singular stratification from the local section $\sigma.$  Then :
\begin{enumerate}
\item[1)] an isometry $\varphi_{[\beta]}$ defined in  Proposition \ref{prop-holonomia-singular} that let $\Omega$ invariant is the identity,
\item[2)] $\Holsing(\sigma)$ is generated by the reflections  in the hypersurfaces of the singular stratification of the local section. 
\end{enumerate}
\end{prop}

\dem  a) Let $p$ a point of $\Omega.$ Since the leaves are compact, $L_{p}$ intercept  $\Omega$  only a finite number of times. Hence, there exists a number $n_{0}$ such that $\varphi_{[\beta]}^{n_{0}}(p)=p.$ Let $K:=\{\varphi_{[\beta]}^{i}(p)\}_{0\leq i<n_{0}}\subset \Omega.$ 

\begin{lema}
There exist only one ball   $B_{r}(x)\supset K$ with minimal radio  $r$  The centre  $x$ of this ball belongs to   $K.$  
\end{lema}
\textbf{Proof of the lemma.} The proof of the lemma is standard, so we sketch the principal steps.

\afir 1: There exists a ball $B_{r}(x)\supset K$ with minimal radio. The centre $x$ belongs to $K.$ 

This follows from the convexity of the balls.

\afir 2: A ball  $B_{r}(x)\supset K$ with minimal radio is unique. 

To check this statement suppose that there exists two balls  $B_{r}(x_{1})$ and  $B_{r}(x_{2})$ that contain  $K$ and have minimal radio $r$. Let $x_{3}$ be the middle point of the segment that joins  $x_{1}$ to $x_{2}.$ Then is possible to find a radio  $\tilde{r}<r,$ such that  $B_{\tilde{r}}(x_{3})\supset (B_{r}(x_{1})\cap B_{r}(x_{2}))$ $\fimdem$

Now we return to the proof of the item a) of the proposition.

 Since $\varphi_{[\beta]}$ let $K$ be invariant,  $K= \varphi_{[\beta]}(K)\subset B_{r}(\varphi_{[\beta]}(x)).$ Since $B_{r}(\varphi_{[\beta]}(x))$ is the ball with the minimal radio that contains $K,$  then  $\varphi_{[\beta]}$ fixes the point $x\in \Omega.$  On the other hand, since the holonomy of regular leaves are trivial, $d_{x}\varphi_{[\beta]}$ is the identity. Since $\varphi_{[\beta]}$ is an isometry, it is the identity.

b) Let  $\varphi_{[\beta]}\in Holsing (\sigma).$ We can compose  $\varphi_{[\beta]}$ with reflections  $R_{i}$'s in the walls of the singular stratification such that   $R_{1}\circ\cdots\circ R_{k}\circ \varphi_{[\beta]}$ let $\Omega$ invariant  and hence, it follows from the item a)  that $R_{1}\circ\cdots \circ R_{k}\circ \varphi_{[\beta]}$ is the identity. We conclude that $\Holsing(\sigma)$ is generated by the reflections  in the hypersurfaces of the singular stratification. $\fimdem$

\begin{cor}
 Let $\F$ be a s.r.f.s. on a complete riemannian manifold $M.$ Suppose that the leaves are compact and that the holonomies of regular leaves are trivial. Let $\tub(L_{q})$ be a tubular neighborhood of a leaf $L_{q},$  $L_{p}$ a regular leaf that belongs to $\tub(L_{q})$ and $\Pi:\tub(L_{q})\rightarrow L_{q}$ the orthogonal projection. Then $L_{p}$ is the total space of a fiber bundle with a projection  $\Pi,$ a basis   $L_{q}$ and a fiber that is diffeomorphic to an isoparametric submanifold of an euclidian space.
\end{cor}
\dem  $\Pi: L_{p}\rightarrow L_{q}$ is a submersion for the foliation is singular. 
 
\afir $\Pi^{-1}(c)=\Sigma_{c}\cap L_{p}$  has only one connected component. 

To check this statement suppose that  $\tilde{L}_{x},\tilde{L}_{y}\subset \Sigma\cap L_{p}$ are two disjoint leaves of  $\Sigma\cap \F.$ 
We can suppose that $x$ and $y$ belong to the same local section.
Since   $x,y\in L_{p},$ there exists  $\varphi_{[\beta]}\in \Holsing(\sigma)$ such that  $\varphi_{[\beta]}(x)=y.$ The above corollary implies that  $\varphi_{[\beta]}$ is a composition of reflections in the hypersurfaces of the singular stratification and hence  $y=\varphi_{[\beta]}(x)\in \tilde{L}_{x}.$ Therefore $\tilde{L}_{x},\tilde{L}_{y}$ are the same leaf. 

 Now  the proposition follows form the slice theorem and a theorem of  Ehresmann \cite{Eh}, which we recall bellow. 

\emph{ Let $\Pi: L\rightarrow K$ a submersion, where $L$ and $K$ are compact manifolds. Suppose that  $\Pi^{-1}(c)$ has only one connected component for each value $c.$ Then the preimages are each other diffeomorphic  and $\Pi:L\rightarrow K$ is the projection of a fiber bundle with total space $L,$ basis $K$ and fiber  $\Pi^{-1}(c).$} $\fimdem$

\begin{prop}
\label{estrutura-transversa-do-fecho} 

Let $\F$ be a s.r.f.s. on a complete riemannian manifold $M,$ $\sigma$ a local section and $p\in\sigma.$ Then
\[\overline{\Holsing(\sigma)}\cdot p= \overline{L_{p}\cap\sigma}.\]
In other words, the closure of $L_{p}\cap\sigma$ is an orbit of 
complete close pseudogroup of local isometries. In particular $\overline{L_{p}\cap\sigma}$ is a closed submanifold.
\end{prop}

\dem This result follows direct from results of E. Salem about pseudogroups of isometries (see appendix D in\cite{Molino}).

One can argue like Salem (see Proposition 2.6 in \cite{Molino}) to prove that $\overline{\Holsing(\sigma)}$ is complete and closed  for the $C^{1}$ topology. It follows from Theorem 3.1 in \cite{Molino} that 
a complete closed pseudogroup of isometry is a Lie pseudogroup. It also follows from E. Salem that a orbit of this Lie pseudogroup is a closed submanifold (see Corollary 3.3 in \cite{Molino}). Therefore $\overline{\Holsing(\sigma)}\cdot p$ is a closed submanifold. Now it is easy to see that $\overline{\Holsing(\sigma)}\cdot p\supset \overline{\Holsing(\sigma)\cdot p}.$ It is also easy to see that $\overline{\Holsing(\sigma)}\cdot p\subset  \overline{\Holsing(\sigma)\cdot p}.$ To finish the proof we have only to remember that  $\Holsing(\sigma)\cdot p= L_{p}\cap \sigma.$ $\fimdem$

%%%%%%%%%%%%%%%
\section{Examples} 

In this section we illustrate some properties of the singular holonomy constructing singular riemannian foliations with singularities by suspension of a homomorphism. 

We start by recalling what a suspension is. For more details see for example the book of Molino \cite{Molino}[page 28,29; 96,97]. 

Let $B$ and $T$ be riemannian manifolds with dimension $p$ and $n$ respectively and let $\rho: \pi_{1}(B, b_{0})\rightarrow  Iso(T)$ be a homomorphism from the fundamental group of $B$ to the group of isometries of $T.$ Let $\hat{P}:\hat{B}\rightarrow B$ be the projection of the universal cover of $B$ into $B.$  
Then we can define an action of $\pi_{1}(B,b_{0})$ on $\widetilde{M}:=\hat{B}\times T$ as
\[ [\alpha]\cdot(\hat{b},t):=([\alpha]\cdot\hat{b},\rho(\alpha^{-1})\cdot t),\]
where $[\alpha]\cdot \hat{b}$ denote  the  deck transformation associated to $[\alpha]$ applied to a point $\hat{b}\in\hat{B}.$

 We denote the set of orbits of this action by $M$ and the canonical projection by $\Pi:\widetilde{M}\rightarrow M.$  It's possible to see that $M$ is a manifold. Indeed, given a simple open neighborhood $U_{j}\subset B,$ we can construct the following bijection:  

\[\begin{array}{lrll}
\Psi_{j}:& \Pi(\hat{P}^{-1}(U_{j})\times T)& \rightarrow & U_{j}\times T \\
             & \Pi(\hat{b},t)                      & \rightarrow & (\hat{P}(\hat{b})\times t).
\end{array}\]

If  $U_{i}\cap U_{j} \neq \emptyset$ and connected, we can see that 
\[\Psi_{i}\cap\Psi_{j}^{-1}(b,t)= (b, \rho([\alpha]^{-1})t)\]
for a fixed  $[\alpha].$ 

So there exists an unique manifold structure on $M$ for which $\Psi_{j}$ are local diffeomorphisms. We define a map $P$ as 
\[ \begin{array}{rlcl}
P: &M &\rightarrow& B\\
   &\Pi(\hat{b},t)& \rightarrow & \hat{P}(\hat{b})
\end{array}
\]
It's possible to see that $M$ is a total space of a fiber bundle,   $P$ is the   projection of this fiber bundle, $T$ is the fiber, $B$ is the basis and the image of $\rho$ is the structure group. 

At last we define $\F:=\{\Pi(\hat{B},t)\},$ i.e., the projection of the trivial foliation defined as the product of  $\hat{B}$ with each $t.$ It is possible to see that this is a foliation transverse to the fibers of the fiber bundle. 
In addition this foliation is a riemannian foliation such that the transversal metric coincide with the metric of $T.$ 

\begin{ex}
\emph{In what follows we construct a singular riemannian foliation with sections such that the intersection of a local section with the closure of a regular leaf is  an orbit of an action of a subgroup of isometries of the local section. This illustrates  Proposition  \ref{estrutura-transversa-do-fecho}. } 
\end{ex}

Let $T$ denote the product   $\mathbf{R}^{2}\times\mathbf{S}^{1}$ and   $\hat{\F}_{0}$ the singular foliation of codimension 2 on  $T$ such that each leaf is the product of a point of $\mathbf{S}^{1}$ with a circle  in $\mathbf{R}^{2}$ whose   centre is $(0,0).$  It is easy to see that the foliation $\hat{\F}_{0}$ is a singular riemannian foliation with sections and that sections are cylinders. Let $B$ be the circle $\mathbf{S}^{1}$ and $q$ be a irrational number. Then we  define the homomorphism $\rho$ as 

\[
\begin{array}{llcl}
\rho:& \pi_{1}(B,b_{0})&\rightarrow& \Iso (T)\\
     &  n               &\rightarrow& ((x,s)\rightarrow (x,\exp(i\,n\,q)\cdot s)).\\
    \end{array}\]

Finally we define $\F:=\Pi(\hat{B}\times \hat{\F_{0}}).$  One can notice that $\F$ is a singular riemannian foliation with sections such that the intersection of each section with the closure of a regular leaf is an orbit of an isometric action on the section.  Indeed one can see this action as translations along the meridians of a cylinder, which  is a section of the foliation.

\begin{ex}
\emph{In what follows we construct a singular riemannian foliation with sections such that $\Holsing(\sigma)$ has an element that can not be  generated by the reflections  in the hypersurfaces of the singular stratification.}
\end{ex}

Let $T$ be a compact Lie group (e.g. $T=SU(3)$) and a manifold $B$ such that  $\pi_{1}(B)=\mathbf{Z}_{2}$ (e.g. $B=SO(n)$). We define the homomorphism $\rho$ as follows 

\[
\begin{array}{llcl}
\rho:& \pi_{1}(B,b_{0})&\rightarrow& \Iso (T)\\
     &  0               &\rightarrow& (t\rightarrow t)\\
     &  1               &\rightarrow& (t\rightarrow t^{-1}).
\end{array}\]

Let us consider the action of $T$ on itself by conjugation, i.e. $t\cdot g:=t\,g\,t^{-1}.$ The orbits of this action are leaves of a singular riemannian foliation that has tori as sections. We denote this singular foliation by $\hat{F}_{0}.$ It's easy to see that $(T\cdot g)^{-1}= T\cdot g^{-1}.$ This assure us that  
$\F:= \Pi(\hat{B}\times\hat{\F}_{0})$ is a singular foliation on  $M.$ We can give a metric to $M$ such that the  metric of the fibers coincide with the metric of $T.$ Then $\F$ turns out to be a singular riemannian foliation whose sections are contained in the fibers. This sections are tori. 

Now it's  possible to see that the leaves of $\F$ intersect a Wely chamber of each torus in more then one point. In fact 
give a point $x_{1}$ belonging to  a Wely chamber, we can reflect it in the walls of the singular stratification and get another point $x_{2}$ belonging to another Wely chamber and such that $x_{2}^{-1}$ belongs to the same Wely chamber of $x_{1}.$ Since inverse points belong to the same leaf, $x_{2}^{-1}$ belongs to the same leaf of $x_{1}.$

%%%%%%%%%%%%%%%
%%%%%%%%%%%%%%%

%%%%%%%%%%%%%%%

\end{document}